\documentclass[10pt]{amsart}

\usepackage{latexsym,amssymb,amsfonts,amscd,amsthm,amstext}
\unitlength=1mm\linethickness{0.8pt}

\def\newpic#1{}
\numberwithin{equation}{section}

\newtheorem{theorem}{Theorem}[section]
\newtheorem{proposition}{Proposition}[section]

\theoremstyle{definition}

\newtheorem{definition}{Definition}[section]
\newtheorem{example}{Example}[section]
\newtheorem{remark}{Remark}[section]

\newcommand{\be}{\begin{equation}}

\newcommand{\bt}{\begin{tabular}{c}}
\newcommand{\et}{\end{tabular}}

\begin{document}

\author{E. Graczy\'{n}ska}
\address{Technical University of Opole, Institute of Mathematics
\newline
ul. Luboszycka 3, 45-036 Opole, Poland} \email{egracz@po.opole.pl
\hspace{1cm} http://www.egracz.po.opole.pl/}
\author{D. Schweigert}
\address{Technische Universit\"{a}t Kaiserslautern, Fachbereich Mathematik
\newline
Postfach 3049 \\ 67653 Kaiserslautern, Germany}
\email{schweige@mathematik.uni-kl.de}

\title{Fluid varieties}

\begin{abstract}
This is a paper presented at the Universal Algebra Conference in Szeged
University (Hungary) on 4-8 July 2005.

{\it Derived varieties} were invented by P. Cohn in the 60ties.
{\it Derived varieties of a given type} invented D. Schweigert in
\cite{10}.

We want to advertize this notion  as a tool for exploring
the lattice of subvarieties of a given variety. A modification of the
notion of a {\it fluid variety} invented in \cite{10} as a counterpart
of the notion  for a {\it solid variety} (see \cite{4}) is given.
A fluid variety $V$ has no proper derived variety as a subvariety.
We examine some properties of derived and fluid varieties in the lattice
of all varieties of a given type $\tau$. Examples of such varieties of bands
are presented.

\emph {Keywords}: fluid  varieties, derived algebras, derived varieties,
proper varieties.    \\

AMS Mathematical Subject Classification 2000: Primary: 08B99, 08A40,

Secondary: 08B05, 08B15.

\end{abstract}

\maketitle
\section{Notations}
\begin{definition}
Let $V$ be a variety of a fixed type $\tau$: $I \rightarrow N$, where $T$ is an index set and
$N$ is the set of all natural numbers. Then $F$ denotes the set of all fundamental operations
$F = \{f_{i}: i \in I\}$ of type $\tau$, i.e. $\tau(i)$ is the arity of the
operation symbol $f_{i}$, for $i \in I$. Let $\sigma = (t_{i}: i \in I)$ be
a fixed choice of terms of $V$ with $\tau(t_{i}) = \tau(f_{i})$, for every $i \in I$.

For any algebra ${\bf A} = (A, \Omega) = (A, (f_{i}^{\bf A}: i \in I)) \in V$,
of type $\tau$, the algebra ${\bf A}_{\sigma}= (A, (t^{\bf A}_{i}: i \in I))$ or shortly
${\bf A}_{\sigma} = (A, \Omega_{\sigma})$, for $\Omega_{\sigma} = (t_{i} :i \in I)$
is called a {\it derived algebra } (of a given type $\tau$) of ${\bf A}$,
corresponding to $\sigma$ (cf. \cite{10}).
\end{definition}
Recall from \cite{4} (cf. \cite{1}), that for a given $\sigma$, the extension
of $\sigma$ to the map from the set $T(\tau)$ of all term symbols of type
$\tau$ to $T(\tau)$, leaving all the variables unchanged and acting on
composed terms as:
\begin{center}
$\sigma(f_{i}(p_{0},...,p_{n-1})) = \sigma(f_{i})(\sigma(p_{0}),...,\sigma(p_{n-1}))$
\end{center}
is called a {\it hypersubstitution} of type $\tau$  and will be
denoted by $\sigma$.

The hypersubstitution $\sigma$ will be called trivial, if it is the identity
mapping.

The set of all hypersubstitutions of type $\tau$ will
be denoted by $H(\tau)$.
\begin{definition}
The variety generated by the class of all derived algebras ${\bf A}_{\sigma}$,
of algebras ${\bf A} \in V$ will be called the {\it derived variety} of $V$
using $\sigma$ and it will be denoted by $V_{\sigma}$.

For a class $K$ of algebras of a given type $\tau$, $D(K)$ denotes the class
of all derived algebras of $K$ for all possible choices of $\sigma$ of type
$\tau$.
\end{definition}
Let us note, that $V_{\sigma} = HSP(\sigma(V))$, for a given variety $V$ and $\sigma$,
where $\sigma(V)$ denotes the class of all derived algebras ${\bf A}_{\sigma}$,
for ${\bf A} \in V$. $D$ is a class operator examined in \cite{4}
(cf. \cite{9}, \cite{10}).

We accept the following definition from \cite{10}:
\begin{definition}
A derived algebra ${\bf A}_{\sigma}$ is called {\it proper} if ${\bf A}_{\sigma}$
is not isomorphic to ${\bf A}$.
\end{definition}
In \cite{10} many examples of {\it proper algebras} were given.
\begin{definition}
A derived variety $V_{\sigma}$ is {\it proper} if $V_{\sigma}$ does not
equal to $V$, i.e.$V_{\sigma} \neq V$.
\end{definition}
Note, that $V_{\sigma}$ may be proper only for nontrivial $\sigma$.

Recall from \cite{4}:
\begin{definition}
A variety $V$ of type $\tau$ is {\it solid} if $V$ contains all derived
varieties $V_{\sigma}$ for every choice of $\sigma$ of type $\tau$.
\end{definition}
\begin{definition}
A variety $V$ of type $\tau$ is {\it fluid} if the variety $V$ contains
no proper derived varieties $V_{\sigma}$ for every choice of $\sigma$ of
type $\tau$.
\end{definition}
Our aim is to show that fluid varieties appear naturally in many well
known examples of varieties of algebras. Moreover, derived varieties
are an important tool for describing the lattice of all subvarieties of
a given variety and therefore we expect many practical applications
of the invented notions.

Note, that our main definition of fluid variety does not coincide with those
of \cite{10}.

Recall from [5] that:
\begin{definition}
A variety $V$ is {\it trivial} if all algebras in $V$ are {\it trivial}
(i.e. one-element). A subclass $W$ of a variety $V$ which is also a variety is
called {\it subvariety} of $V$. $V$ is {\it minimal} (or {\it equationally complete})
variety if $V$ is not trivial but the only subvariety of $V$, which is not
equal to $V$ is trivial.
\end{definition}

\begin{theorem}
Minimal varieties are fluid.
\end{theorem}
{\it Proof}. Let $V$ be a minimal variety. Assume that
$V_{\sigma} \subseteq V$ for some  derived variety $V_{\sigma}$.
As $V$ is nontrivial, therefore $V_{\sigma}$ is nontrivial and we conclude
that $V_{\sigma} = V$ and therefore $V_{\sigma}$ is not proper. Therefore
$V$ is fluid.
\begin{theorem}
Each nontrivial variety contains a  nontrivial fluid variety.
\end{theorem}
This follows immediately from the fact that each nontrivial variety
contains a minimal variety.

From the definitions above it follows that every trivial variety $0$
is fluid and solid. Similarly, the variety $1$ defined by an empty
set of identities of a given type $\tau$ is also fluid and solid.
Our last example will show that not only trivial varieties may be
fluid and solid.

Natural examples of nontrivial fluid varieties arises:
\begin{example}\label{E:1.1}
The variety of all semilattices in type (2) is fluid.
\end{example}
\begin{example}\label{E:1.2}
The variety of distributive lattices in type (2,2) is fluid.
\end{example}
\begin{example}\label{E:1.3}
The variety of Boolean algebras in type (2,2,1,0,0) is fluid.
\end{example}

First we concentrate on varieties of {\it bands} (cf. \cite{2},
\cite{fen1}, \cite{fen2}, \cite{3}) defined as the variety of all
idempotent semigroups of type (2). The reason is that:

First we point out a practical theorem:
\begin{theorem}
Given algebra ${\bf A}$ and the hypersubstitution $\sigma$ of type $\tau$.
Then an identity $p = q$ of type $\tau$ is satisfied in the derived
algebra ${\bf A}_{\sigma}$ if and only if the derived identity
$\sigma(p) = \sigma(q)$ holds in ${\bf A}$.
\end{theorem}
{\it Proof}. To show this let us notice, that the realization of a term
$p$ of type $\tau$ in the derived algebra ${\bf A}_{\sigma}$ equals to
$\sigma(p)$. This fact can be easily proved by induction on the complexity
of term p.

From the theorem above, it  immediately follows:
\begin{theorem}
Let  $V$ be a variety and given  two hypersubstitutions $\sigma_{1}$ and
$\sigma_{2}$ of type $\tau$.
If $\sigma_{1}(f_{i}) = \sigma_{2}(f_{i})$, is an identity of $V$ for
every $i \in I$, then the derived varieties $V_{\sigma_{1}}$ and
$V_{\sigma_{2}}$ are equal.
\end{theorem}
{\it Proof}. From the assumption it follows that $\sigma_{1}(p) =
\sigma_{2}(p)$, is an identity of $V$ for every polynomial symbol $p$
of type $\tau$. Therefore any identity $p = q$ is satisfied in
${\bf A}_{\sigma_{1}}$ if and only if it is satisfied in ${\bf
A}_{\sigma_{2}}$, for every algebra ${\bf A}$ of $V$.

\begin{definition}
Bands is a variety $B$ of algebras of type (2), defined by: associativity and
idempotency (i.e. {\bf band} is an idempotent semigroup).
\end{definition}

Recall from \cite{1}, p.11:

\begin{proposition}
A derived variety of bands is may not be a variety of bands.
\end{proposition}

\begin{theorem}
The variety of $B$ all bands constitute a not fluid and not
solid variety of type (2).
\end{theorem}
{\it Proof}. For a  $\sigma$, generated by the first (second) projection,
the derived variety $\sigma(B)$ is the variety of left (or right) zero
semigroups. Therefore $B$ is not fluid.

From the proof of the theorem above we conclude the following:
\begin{theorem}
A variety $V$ of bands is fluid if and only if it is minimal.
\end{theorem}

\begin{remark}
There are solid varieties, which are not fluid.
\end{remark}

\begin{proposition}
The variety $W_{1}$ of rectangular bands is solid and not fluid.
\end{proposition}
{\it Proof}. Recall  (cf. \cite{fen1}, \cite{fen2}) that $W_{1}$ is
defined in the variety of $B$ bands by the identity: $y = yxy$. The fact
that this variety is solid was proved in \cite{1}, p. 96. To show
that it is not fluid, given the hypersubstitutions $\sigma_{1}$ and
$\sigma_{2}$ defined by the first and second projections, respectively,
i.e. $\sigma_{1}(xy) = x$ and $\sigma_{2}(xy) = y$.  Take $\sigma =
\sigma_{1}$. Let ${\bf A}_{\sigma} \in \sigma(W_{1})$, for an
algebra ${\bf A} \in W_{1}$. Then the identity $xy = y$ is satisfied
in ${\bf A}_{\sigma}$, as the identity: $\sigma(xy) = y = y =
\sigma(y)$ is satisfied in ${\bf A}$. Therefore the derived variety
of bands $\sigma(W_{1})$ is defined by the identity $xy = y$. We got
that: $\sigma(W_{1})$ is proper derived variety with $\sigma(W_{1})
\subseteq W_{1}$ and therefore $W_{1}$ is not fluid.

\begin{remark}
There are many not solid varieties, which are not fluid.
\end{remark}
\begin{proposition}
The varieties $V_{1}$ and $V_{2}$ of bands defined by the identities:

(1)  $zxy=zyx$ \hspace{5mm} and \hspace{5mm} (2) $yxz = xyz$, respectively,

are mutually derived and not fluid. Moreover, they are not solid.
\end{proposition}

{\it Proof}. This follows from the fact, that they are mutually derived via the
hypersubstitution $\sigma$ generated by $\sigma(xy) = yx$, i.e.
$\sigma(V_{1}) = V_{2}$ and $V_{2} = \sigma(V_{1})$. To show that
$\sigma(V_{1}) = V_{2}$, given an algebra ${\bf A}_{\sigma} \in \sigma(V_{1})$,
then the identity $yxz = xyz$
is satisfied in ${\bf A}_{\sigma}$ as: $\sigma(xyz) = \sigma(z)\sigma(xy) =
zyx = \sigma(yxz) = \sigma(z)\sigma(yx) = zxy$ is satisfied in $V_{1}$ and
vice versa.
Moreover, consider $\sigma$ given by: $\sigma(x,y) = x$. We show that then
the variety $\sigma(V_{1})$ is the variety defined by $yx = y$. To show
this given an algebra ${\bf A}_{\sigma}$, for ${\bf A} \in V_{1}$. Then:
$\sigma(yx) = y = \sigma(y)$  is satisfied in ${\bf A}$. We got that
$\sigma(V_{1}$) is different as $V_{1}$ and therefore $\sigma(V_{1})$ is
proper and $\sigma(V_{1}) \subseteq V_{1}$ (cf. \cite{fen1}, \cite{fen2})
which proves  that $V_{1}$ is not fluid. Similarly $V_{2}$ is not fluid.
Both of them are not solid, as the second and first projections, for
$V_{1}$ and $V_{2}$ respectively, give rise to a trivial identity $x = y$.

The next three propositions show some regularities in the diagram describing
the lattice of all identities of bands in \cite{fen1}, p. 244:

\begin{proposition}
The varieties $V_{3}$ and $V_{4}$ of bands defined by the identities:

(3)  $yx = yxy$ \hspace{5mm} and \hspace{5mm} (4) $xy = yxy$, respectively,

are mutually derived and not fluid.  Moreover, they are not solid.
\end{proposition}
{\it Proof}. By the similar argument as in the previous example, let us note that
the varieties $V_{3}$ and $V_{4}$ are mutually derived via the
hypersubstitution $\sigma$ given by $\sigma(xy) = yx$, i.e.
$\sigma(V_{3}) = V_{4}$ and $\sigma(V_{4}) = V_{3}$. To show the first one
equation, given algebra ${\bf A}_{\sigma} \in \sigma(V_{3})$, for an algebra
${\bf A} \in V_{3}$. Then the identity $xy = yxy$ is satisfied in
${\bf A}_{\sigma}$, as the identity:
$\sigma(xy) = yx = yxy = \sigma(yxy)$ is satisfied in $V_{3}$ and therefore
in ${\bf A}$. Similarly one can show that $\sigma(V_{4}) = V_{3}$.

Let now $\sigma$ be given by the first projection, i.e. $\sigma(xy) = x$.
Oncemore, the variety $\sigma(V_{3})$ is the variety of bands defined by
$yx = y$, as given $\sigma({\bf A}) \in \sigma(V_{3})$, for ${\bf A} \in V_{3}$.
Then the identity $yx = y$ is satisfied in ${\bf A}_{\sigma}$, as:
$\sigma(yx) = y = y = \sigma(y)$ is satisfied in ${\bf A}$. We got:
$\sigma(V_{3})$ different as $V_{3}$, therefore $\sigma(V_{3})$ is proper and
$\sigma(V_{3}) \subseteq V_{3}$ and $V_{3}$ is not fluid.

Similarly one can see that $V_{4}$ is not fluid. Both are not solid. The
argument is similar as in previous proposistion.

\begin{proposition}
The varieties $V_{5}$ and $V_{6}$ of bands defined by the identities:

(5) $xzy = zxzy$ \hspace{5mm} and \hspace{5mm} (6) $yxz = yzxz$,
respectively,

are mutually derived and not fluid. Moreover, they are not solid.
\end{proposition}
{\it Proof}. Similarly as in the previous proposition, $V_{5} =
\sigma(V_{6})$ and $V_{6} = \sigma(V_{5})$ for $\sigma(xy) = yx$.
Therefore $V_{5}$ and $V_{6}$ are mutually derived. Both are not
fluid, as they contain the proper derived varieties:
$\sigma_{1},_{2}(V_{5})$ or $\sigma_{1},_{2}(V{6})$, respectively,
for $\sigma_{1}$, $\sigma_{2}$ being the first or the second
projection.

\begin{proposition}
The variety $W_{2}$ defined by the identity: (7) $zxyz = zyxz$ is solid
and not fluid.
\end{proposition}
{\it Proof}. For solidity confront \cite{1}, p. 96. As $W_{2}$ is not minimal,
therefore we conclude that it is not fluid.

In fact:

\begin{theorem}
A solid variety $W$ is fluid if and only if $W = W_{\sigma}$, for
every $\sigma \in H(\tau)$.
\end{theorem}
{\it Proof}. If $W$ is solid, then every derived variety
$W_{\sigma}$ is included in $W$, as $W$ contains all derived
algebras of a given type. Therefore $W$ is fluid if and only if all
its derived varieties are not proper.

We shall express the situation of Propostions  1.3 - 1.6 in the
diagram, which describes the bottom part of the lattice of all
identities of bands, see \cite{fen1} and \cite{fen2} p. 244 and
confront Propostion 3.1.5 of \cite{1}, p. 77:

\vspace{5mm}

\input{fluid1.pic}

\end{document}